\documentclass[graybox]{svmult}
\usepackage{mathptmx}       
\usepackage{helvet}         
\usepackage{courier}        
\usepackage{type1cm}        
\usepackage{makeidx}         
\usepackage{graphicx}        
\usepackage{multicol}        
\usepackage[bottom]{footmisc}
\usepackage{amsmath, amsfonts, amssymb}
\usepackage{color}
\usepackage{tikz}

\makeindex             
                       
\begin{document}

\title*{A Hybrid High-Order method for the convective Cahn--Hilliard problem in mixed form}
\titlerunning{A HHO method for the Cahn--Hilliard problem}
\author{Florent Chave, Daniele A. Di Pietro and Fabien Marche}
\institute{
Florent Chave, Daniele A. Di Pietro and Fabien Marche
\at
University of Montpellier, Institut Montpelli\'{e}rain Alexander Grothendieck\\
34095 Montpellier, France\\
\email{\{florent.chave,~daniele.di-pietro,~fabien.marche\}@umontpellier.fr}
}
\maketitle

\abstract{We propose a novel Hybrid High-Order method for the Cahn--Hilliard problem with convection. The proposed method is valid in two and three space dimensions, and it supports arbitrary approximation orders on general meshes containing polyhedral elements and nonmatching interfaces. An extensive numerical validation is presented, which shows robustness with respect to the P\'eclet number.}
\keywords{Hybrid High-Order, Cahn--Hilliard equation, phase separation, mixed formulation, polyhedral meshes, arbitrary order
\\[5pt]
{\bf MSC }(2010){\bf:} 65N08, 65N30 65N12}

\section{Cahn--Hilliard equation}

Let $\Omega\subset\mathbb{R}^d$, $d\in\{2,3\}$, denote a bounded connected convex polyhedral domain with Lipschitz boundary $\partial\Omega$ and outward normal $\vec{n}$, and let $t_F>0$. The convective Cahn--Hilliard problem consists in finding the order-parameter $c : \Omega \times \left(0,t_F\right] \to \mathbb{R}$ and the chemical potential $w : \Omega \times \left(0,t_F\right] \to \mathbb{R}$ such that
\begin{subequations}\label{eq:cont}
  \begin{alignat}{2}\label{eq:cont:1}
    {\rm d}_t c - \frac{1}{{\rm Pe}}\Delta w + \nabla\cdot (\vec{u} c) &= 0 &\qquad&\text{in $\Omega\times(0,t_F\rbrack$}
    \\
    w &= \Phi'(c) - \gamma^2\Delta c &\qquad&\text{in $\Omega\times(0,t_F\rbrack$}
    \\
    c(0) &= c_0 &\qquad&\text{in $\Omega$}
    \\
    \partial_{\vec{n}}c=\partial_{\vec{n}}w&=0 &\qquad&\text{on $\partial\Omega\times (0,t_F\rbrack$}
  \end{alignat}
\end{subequations}
where $\gamma > 0$ is the interface parameter (usually taking small values), Pe$>0$ is the P\'eclet number, $\vec{u}$ the velocity field such that $\nabla\cdot\vec{u}=0$ in $\Omega$ and $\Phi$ the free-energy such that $\Phi(c) := \frac14 (1-c^2)^2$. This formulation is an extension of the Cahn--Hilliard model originally introduced in \cite{Cahn.Hilliard:58,Cahn:61} and a first step towards coupling with the Navier--Stokes equations. 

In this work we extend the HHO method of \cite{Chave.Di-Pietro.ea:16} to incorporate the convective term in (\ref{eq:cont:1}). Therein, a full stability and convergence analysis was carried out for the non-convective case, leading to optimal estimates in $(h^{k+1}+\tau)$ (with $h$ denoting the meshsize and $\tau$ the time step) for the the $C^0(H^1)$-error on the order-parameter and $L^2(H^1)$-error on the chemical potential. The convective term is treated in the spirit of \cite{Di-Pietro.Droniou.ea:14}, where a HHO method fully robust with respect to the P\'eclet number was presented for a locally degenerate diffusion-advection-reaction problem.

The proposed method offers various assets: (i) fairly general meshes are supported including polyhedral elements and nonmatching interfaces; (ii) arbitrary polynomial orders, including the case $k=0$, can be considered; (iii) when using a first-order (Newton-like) algorithm to solve the resulting system of nonlinear algebraic equations, element-based unknowns can be statically condensed at each iteration. 

The rest of this paper is organized as follows: in Section \ref{sec:HHO}, we recall discrete setting including notations and assumptions on meshes, define localy discrete operators and state the discrete formulation of (\ref{eq:cont}). In Section \ref{sec:numerical}, we provide an extensive numerical validation.

\section{The Hybrid High-Order method\label{sec:HHO}}

In this section we recall some assumptions on the mesh, introduce the notation, and state the HHO discretization.

\subsection{Discrete setting}

We consider sequences of refined meshes that are regular in the sense of \cite[Chapter~1]{Di-Pietro.Ern:12}. Each mesh $\mathcal{T}_h$ in the sequence is a finite collection $\{T\}$ of nonempty, disjoint, polyhedral elements such that $\overline{\Omega} = \bigcup_{T\in\mathcal{T}_h} \overline{T}$ and $h = \max_{T\in\mathcal{T}_h} h_T$ (with $h_T$ the diameter of $T$). For all $T\in\mathcal{T}_h$, the boundary of $T$ is decomposed into planar faces collected in the set $\mathcal{F}_T$. For admissible mesh sequences, card($\mathcal{F}_T$) is bounded uniformly in $h$. Interfaces are collected in the set $\mathcal{F}_h^{\text{i}}$, boundary faces in $\mathcal{F}_h^{\text{b}}$ and we define $\mathcal{F}_h := \mathcal{F}_h^{\text{i}} \cup \mathcal{F}_h^{\text{b}}$. For all $T\in\mathcal{T}_h$ and all $F\in\mathcal{F}_T$, the diameter of $F$ is denoted by $h_F$ and the unit normal to $F$ pointing out of $T$ is denoted by $\vec{n}_{TF}$.

To discretize in time, we consider for sake of simplicity a uniform partition $(t^n)_{0\leq n \leq N}$ of the time interval $[0,t_F]$ with $t^0=0$, $t^N = t_F$ and $t^n-t^{n-1} = \tau$ for all $1\leq n\leq N$. For any sufficiently regular function of time $\varphi$ taking values in a vector space $V$, we denote by $\varphi^n\in V$ its value at discrete time $t^n$, and we introduce the backward differencing operator $\delta_t$ such that, for all $1\leq n\leq N$,
\begin{equation*}
\delta_t\varphi^n := \dfrac{\varphi^n - \varphi^{n-1}}{\tau} \in V.
\end{equation*}

\subsection{Local space of degrees of freedom}

For any integer $l\geq 0$ and $X$ a mesh element or face, we denote by $\mathbb{P}^{l}(X)$ the space spanned by the restrictions to $X$ of $d$-variate polynomials of order $l$. Let
$$\underline{U}_h^k:=\Bigg(\times_{T\in\mathcal{T}_h}\mathbb{P}^{k+1}(T)\Bigg)\times\Bigg(\times_{F\in\mathcal{F}_h}\mathbb{P}^k(F)\Bigg)$$
be the global degrees of freedoms (DOFs) space with single-valued interface unknowns. We denote by $\underline{v}_h = ((v_T)_{T\in\mathcal{T}_h},(v_F)_{F\in\mathcal{F}_h})$ a generic element of $\underline{U}_h^k$ and by $v_h$ the piecewise polynomial function such that $v_h|_T = v_T$ for all $T\in\mathcal{T}_h$. For any $T\in\mathcal{T}_h$, we denote by $\underline{U}_T^k$ and $\underline{v}_T = (v_T, (v_F)_{F\in\mathcal{F}_T})$ the restrictions to $T$ of $\underline{U}_h^k$ and $\underline{v}_h$, respectively.
  \begin{figure}\center
    \newdimen\R
    \R=1cm
    \newdimen\r
    \r=0.30cm
    \begin{tikzpicture}[scale=1]
      \begin{scope}[yshift=-3\R]
        \draw (0:\R) \foreach \x in {60,120,...,360} {
          -- node[rotate=\x-120]{$\bullet$} (\x:\R)
        } -- cycle (90:\R) node[above] {$k=0$};
        \foreach \x in {0,120,...,240} {
          \node[anchor=center,color=black!40] at (\x:\r) {$\bullet$};
        };
      \end{scope}
      \begin{scope}[xshift=3\R,yshift=-3\R]
        \draw (0:\R) \foreach \x in {60,120,...,360} {
          -- node[rotate=\x-120]{$\bullet\bullet$} (\x:\R)
        } -- cycle (90:\R) node[above] {$k=1$};
        \foreach \x in {0,60,...,300} {
          \node[anchor=center,color=black!40] at (\x:\r) {$\bullet$};
        };
      \end{scope}
      \begin{scope}[xshift=6\R,yshift=-3\R]
        \draw (0:\R) \foreach \x in {60,120,...,360} {
          -- node[rotate=\x-120]{$\bullet\bullet\bullet$} (\x:\R)
        } -- cycle (90:\R) node[above] {$k=2$};
        \foreach \x in {0,36,...,360} {
          \node[anchor=center,color=black!40] at (\x:\r) {$\bullet$};
        };
      \end{scope}
    \end{tikzpicture}
    \caption{Local DOF space for $k=0,1,2$. Internal DOFs (in gray) can be statically condensed at each Newton iteration.\label{fig:dofs:local}}
  \end{figure}
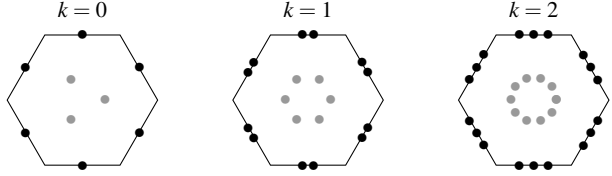%
\subsection{Local diffusive contribution}

Consider a mesh element $T\in\mathcal{T}_h$. We define the local potential reconstruction $\mathrm{p}_T^{k+1}:\underline{U}_T^k \rightarrow \mathbb{P}^{k+1}(T)$ such that, for all $\underline{v}_T := (v_T, (v_F)_{F\in\mathcal{F}_T})\in\underline{U}_T^k$ and all $z\in\mathbb{P}_T^{k+1}$,
\begin{eqnarray*}
(\nabla \mathrm{p}_T^{k+1}\underline{v}_T, \nabla z)_T &=& -(v_T,\Delta z)_T + \sum_{F\in\mathcal{F}_T} (v_F, \nabla z\cdot\vec{n}_{TF})_F,%
\end{eqnarray*}
with closure condition $\int_T (\mathrm{p}_T^{k+1}\underline{v}_T - v_T) = 0$.
We introduce the local diffusive bilinear form $a_T$ on $\underline{U}_T^k\times\underline{U}_T^k$ such that, for all $(\underline{u}_T,\underline{v}_T)\in\underline{U}_T^k\times\underline{U}_T^k$ 
\begin{equation*}
a_T(\underline{u}_T,\underline{v}_T) := (\nabla \mathrm{p}_T^{k+1}\underline{u}_T, \nabla \mathrm{p}_T^{k+1}\underline{v}_T)_T + s_T(\underline{u}_T,\underline{v}_T),
\end{equation*}
with stabilization bilinear form $s_T:\underline{U}_T^k\times\underline{U}_T^k\rightarrow \mathbb{R}$ such that
\begin{equation*}
s_T(\underline{u}_T,\underline{v}_T) := \sum_{F\in\mathcal{F}_T} h_F^{-1}(\pi_F^k(u_F-u_T), \pi_F^k(v_F-v_T))_F,
\end{equation*}
where, for all $F\in\mathcal{F}_h$, $\pi_F^k:L^1(F)\rightarrow \mathbb{P}^ k(F)$ denotes the $L^2$-orthogonal projector onto $\mathbb{P}^ k(F)$.

\subsection{Local convective contribution}

For any mesh element $T\in\mathcal{T}_h$, we define the local convective derivative reconstruction $\mathrm{G}_{\vec{u},T}^{k+1} : \underline{U}_T^k \rightarrow \mathbb{P}^{k+1}(T)$ such that, for all $\underline{v}_T := (v_T, (v_F)_{F\in\mathcal{F}_T})\in\underline{U}_T^k$ and all $w\in\mathbb{P}^{k+1}(T)$,
\begin{eqnarray*}
(\mathrm{G}_{\vec{u},T}^{k+1}\underline{v}_T, w)_T &=& -(v_T,\vec{u}\cdot\nabla w)_T + \sum_{F\in\mathcal{F}_T}(v_F,(\vec{u}\cdot\vec{n}_{TF})w)_F.%
\end{eqnarray*}
The local convective contribution $b_{\vec{u},T}$ on $\underline{U}_T^k\times\underline{U}_T^k$ is such that, for all $(\underline{u}_T,\underline{v}_T)\in\underline{U}_T^k\times\underline{U}_T^k$ 
\begin{equation*}
b_{\vec{u},T}(\underline{u}_T,\underline{v}_T) := (\mathrm{G}_{\vec{u},T}^{k+1}\underline{u}_T, v_T)_T + s_{\vec{u},T}(\underline{u}_T,\underline{v}_T).
\end{equation*}
with local upwind stabilization bilinear form $s_{\vec{u},T}:\underline{U}_T^k\times\underline{U}_T^k\rightarrow \mathbb{R}$ such that  
\begin{equation*}
s_{\vec{u},T}(\underline{u}_T,\underline{v}_T) := \sum_{F\in\mathcal{F}_T} (\frac{|\vec{u}\cdot \vec{n}_{TF}|-\vec{u}\cdot \vec{n}_{TF}}{2}(u_F-u_T), v_F - v_T)_F.
\end{equation*}
Notice that the actual computation of $\mathrm{G}_{\vec{u},T}^{k+1}$ is not required, as one can simply use its definition to expand the cell-based term in the bilinear form $b_{\vec{u},T}$.

\subsection{Discrete problem}

Denote by $\underline{U}_{h,0}^k := \{\underline{v}_h = ((v_T)_{T\in\mathcal{T}_h},(v_F)_{F\in\mathcal{F}_h})\in \underline{U}_h^k | \int_{\Omega} v_h = 0\}$ the zero-average DOFs subspace of $\underline{U}_{h}^k$. We define the global bilinear forms $a_h$ and $b_{\vec{u},h}$ on $\underline{U}_h^k\times\underline{U}_h^k$ such that, for all $(\underline{u}_h,\underline{v}_h)\in\underline{U}_h^k\times\underline{U}_h^k$
\begin{equation*}
a_h(\underline{u}_h,\underline{v}_h) := \sum_{T\in\mathcal{T}_h} a_T(\underline{u}_T,\underline{v}_T),\qquad b_{\vec{u},h}(\underline{u}_h,\underline{v}_h) := \sum_{T\in\mathcal{T}_h} b_{\vec{u},T}(\underline{u}_T,\underline{v}_T).
\end{equation*}
The discrete problem reads: For all $1\leq n \leq N$, find $(\underline{c}_h^n,\underline{w}_h^n)\in\underline{U}_{h,0}^k\times\underline{U}_h^k$ such that
\begin{subequations}
  \begin{alignat*}{2}
    &(\delta_t c_h^n,\varphi_h) + \frac{1}{{\rm Pe}}a_h(\underline{w}^n_h,\underline{\varphi}_h) + b_{\vec{u},h}(\underline{c}^n_h,\underline{\varphi}_h) = 0 &\qquad&\forall\underline{\varphi}_h\in\underline{U}_h^k
    \\
    &(w_h^n,\psi_h) = (\Phi'(c_h^n),\psi_h) + \gamma^2 a_h(\underline{c}^n_h,\underline{\psi}_h) &\qquad&\forall\underline{\psi}_h\in\underline{U}_h^k
  \end{alignat*}
\end{subequations}
where $\underline{c}_h^0\in\underline{U}_{h,0}^k$ solves $a_h(\underline{c}_h^0,\underline{\varphi}_h) = -(\Delta c_0,\varphi_h)$ for all $\underline{\varphi}_h\in\underline{U}_h^k$.

\section{Numerical test cases\label{sec:numerical}}

In this section, we numerically validate the HHO method.

\subsection{Disturbance of the steady solution\label{subsec:steady}}

\begin{figure}[t]\center
\includegraphics[width=3.5cm]{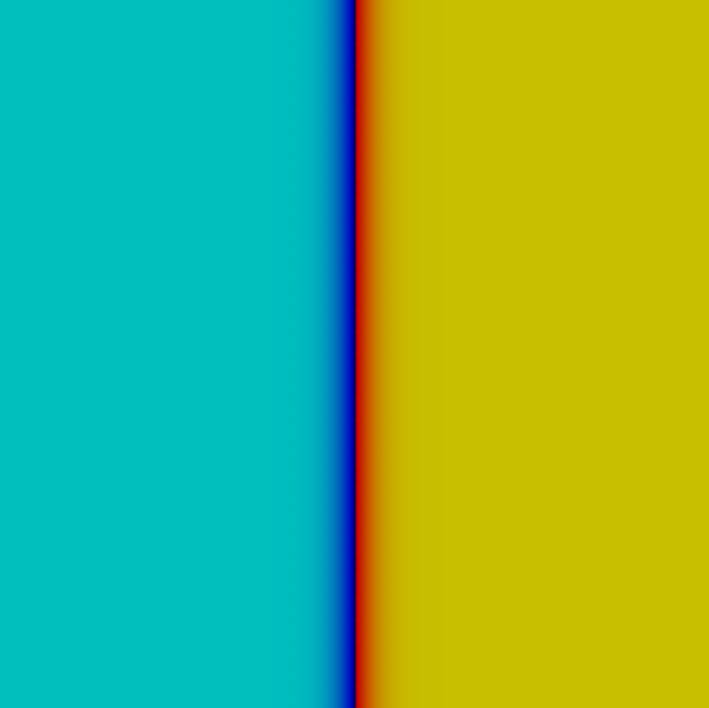}\hspace*{0.5cm}
\includegraphics[width=3.5cm]{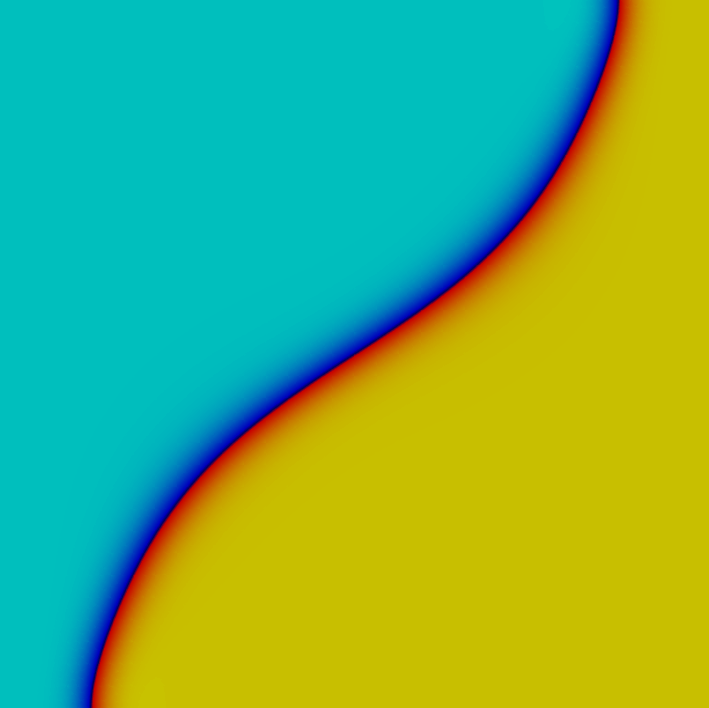}\hspace*{0.5cm}
\includegraphics[width=3.5cm]{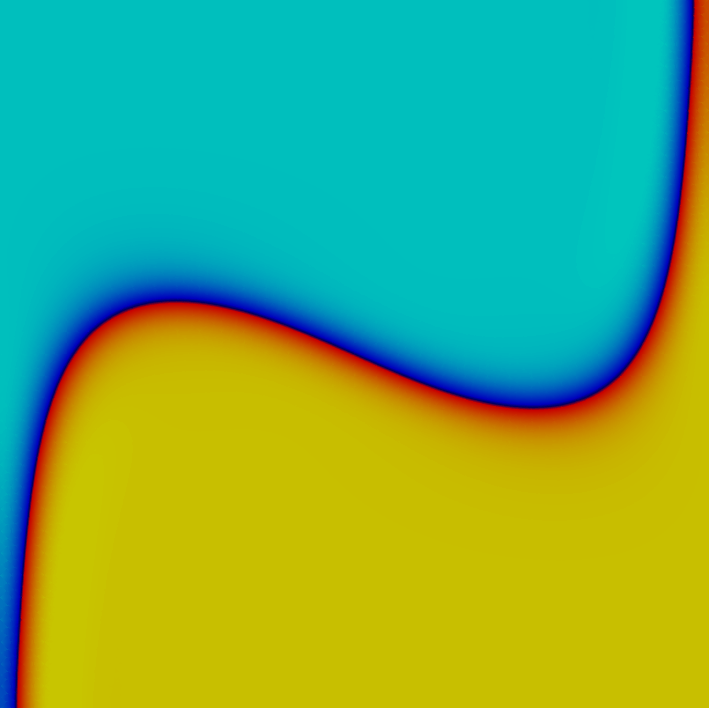}\\[0.5cm]
\includegraphics[width=3.5cm]{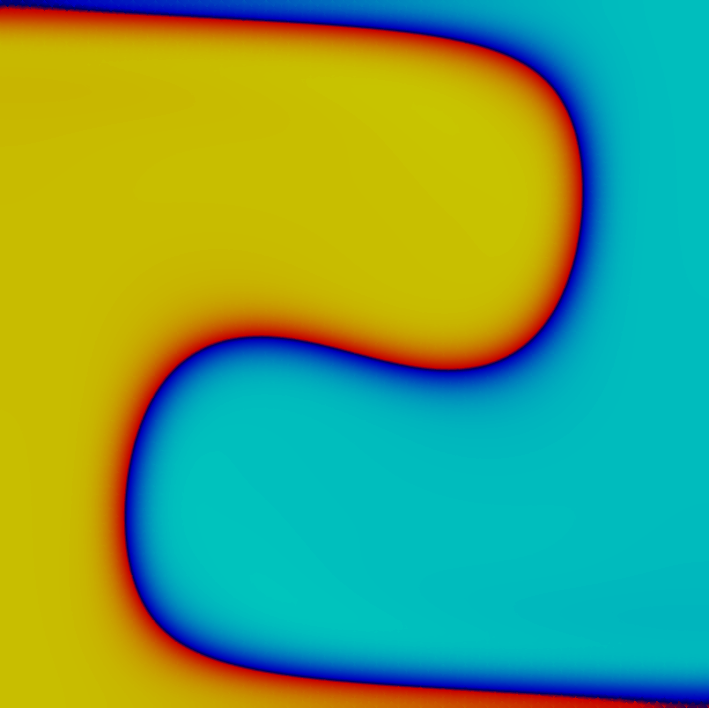}\hspace*{0.5cm}
\includegraphics[width=3.5cm]{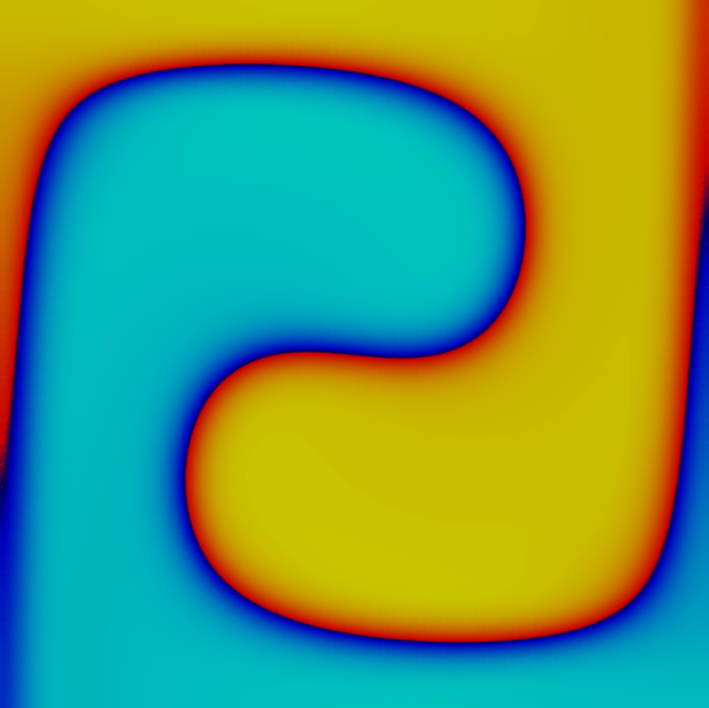}\hspace*{0.5cm}
\includegraphics[width=3.5cm]{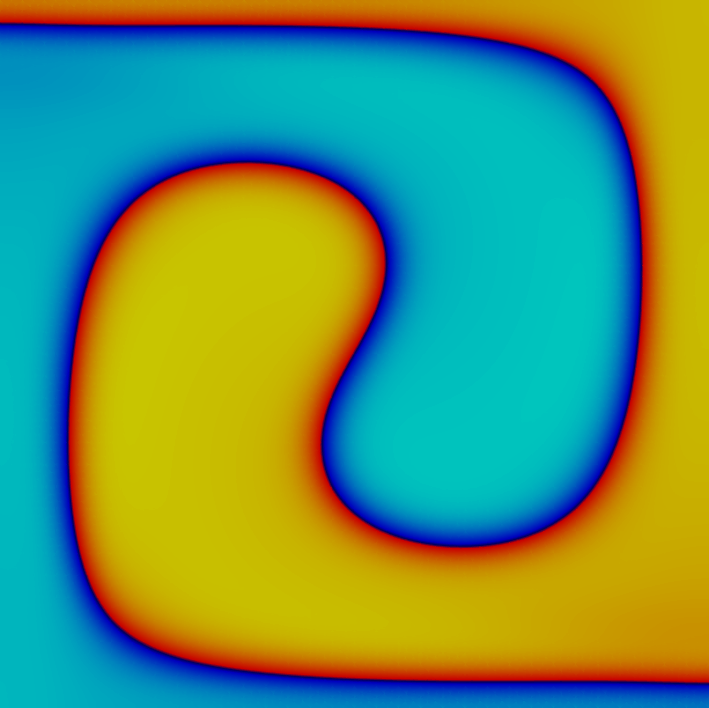}\vfill
\caption{Steady solution perturbed by a circular velocity field (left to right, top to bottom).\label{fig:steady}}
\end{figure}

For the first test case, we use a piecewise constant approximation ($k=0$), discretize the domain $\Omega = (0,1)^2$ by a triangular mesh ($h=1.92\cdot 10^{-3}$) with $\gamma = 5\cdot 10^{-2}$, $\tau = \gamma^2$ and $\mathrm{Pe} = 1$. The initial condition for the order-parameter and the velocity field are given by
\begin{equation*}
c_0(\vec{x}) := \tanh (\dfrac{2x_1 - 1}{2\sqrt{2}\gamma^2}), \quad \vec{u}(\vec{x}) := 20\cdot\begin{pmatrix}x_1(x_1-1)(2x_2-1)\\-x_2(x_2-1)(2x_1-1) \end{pmatrix},\quad \forall\vec{x}\in\Omega.
\end{equation*}
The result is depicted in Figure \ref{fig:steady} and shows that the method is well-suited to capture the interface dynamics subject to a strong velocity fields. 


\subsection{Thin interface between phases\label{subsec:thin}}

\begin{figure}[t]\center
\includegraphics[width=3.5cm]{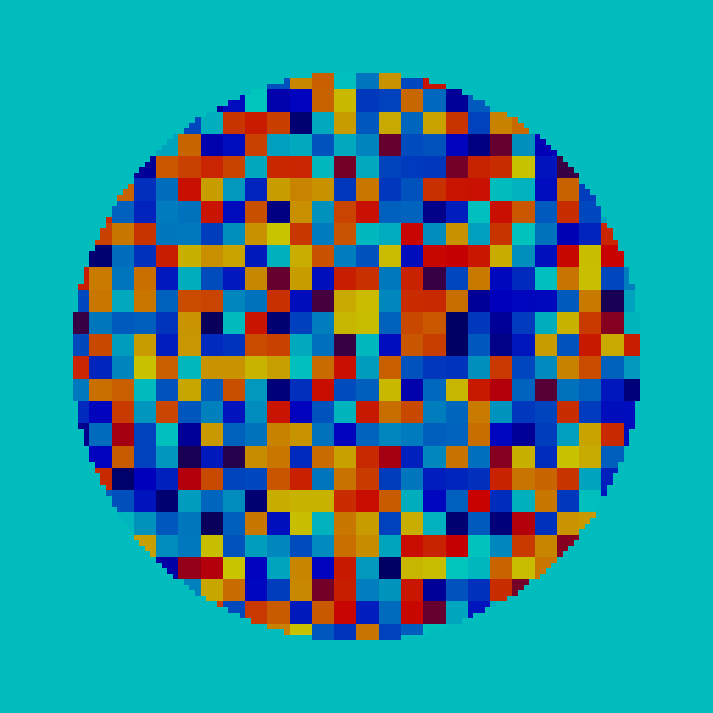}\hspace*{0.5cm}
\includegraphics[width=3.5cm]{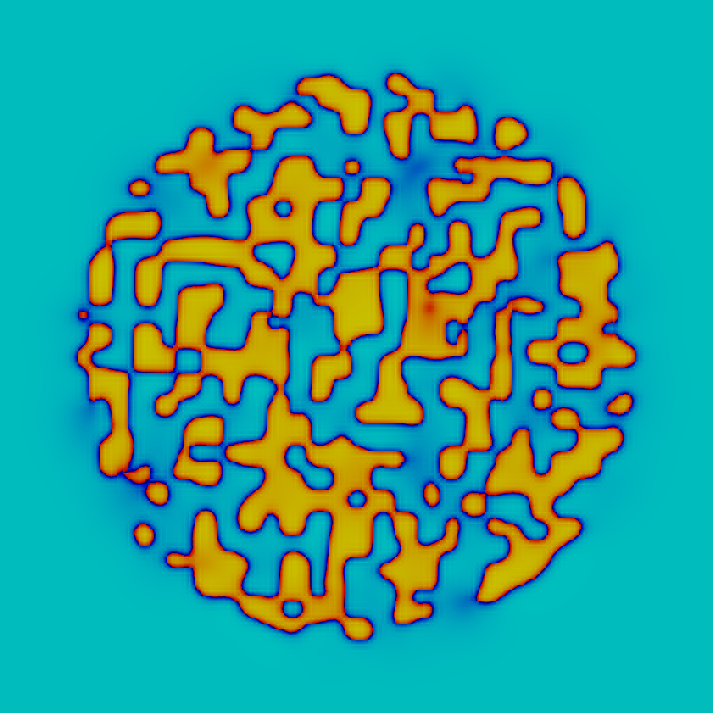}\hspace*{0.5cm}
\includegraphics[width=3.5cm]{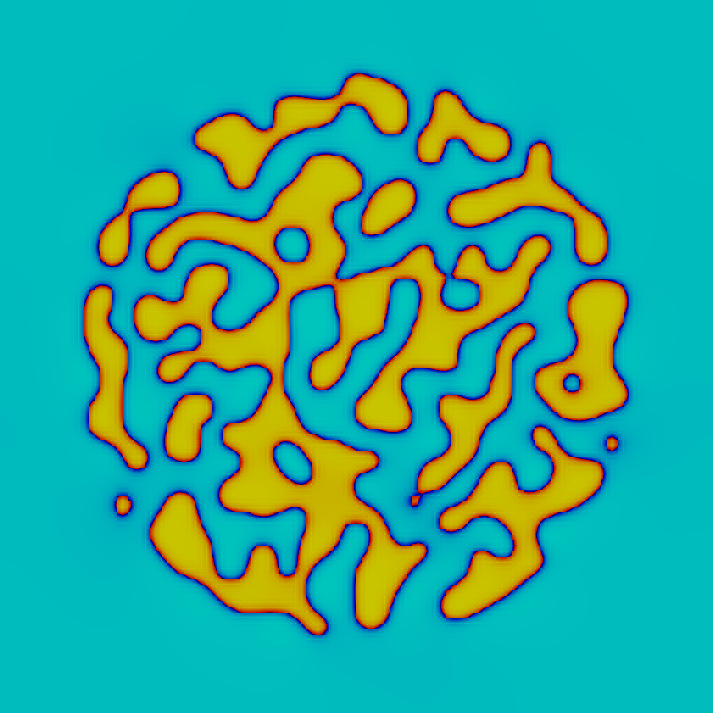}\\[0.5cm]
\includegraphics[width=3.5cm]{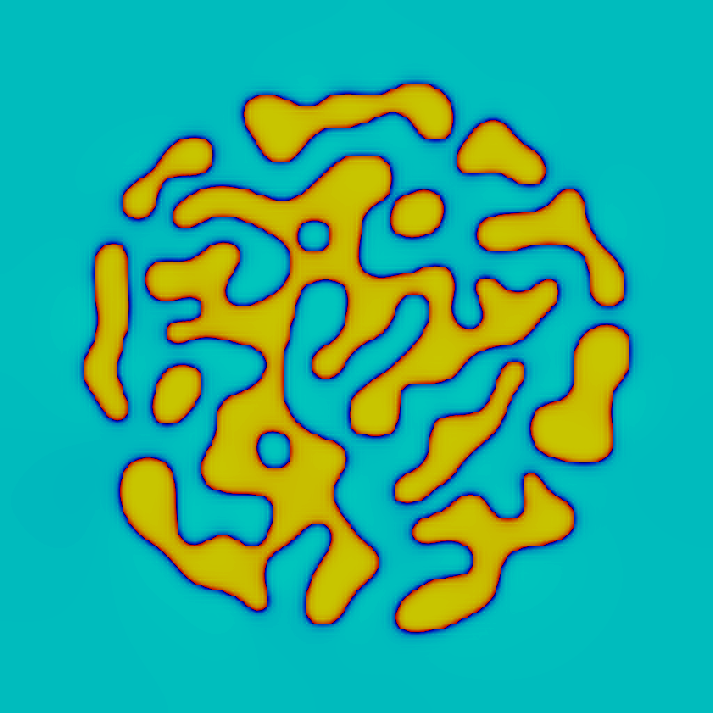}\hspace*{0.5cm}
\includegraphics[width=3.5cm]{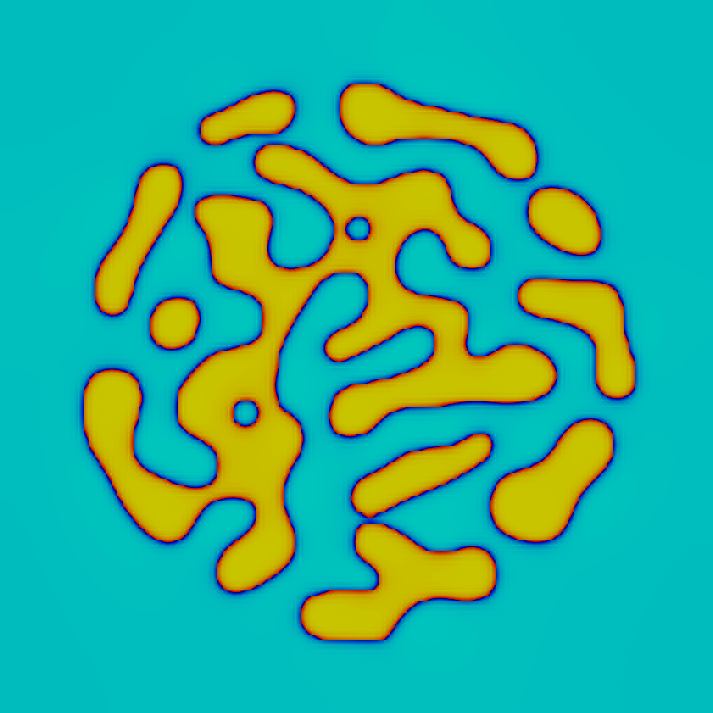}\hspace*{0.5cm}
\includegraphics[width=3.5cm]{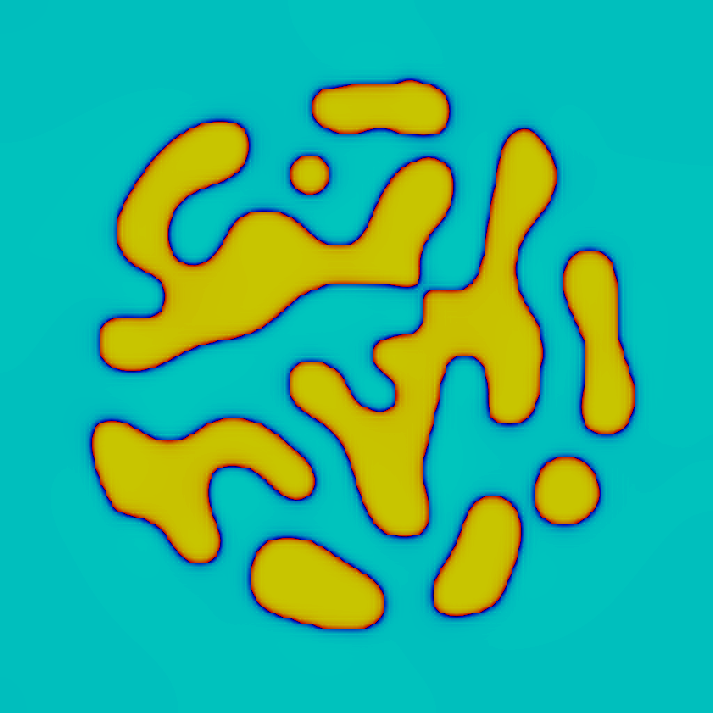}
\caption{Evolution of spinodal decomposition with thin interface (left to right, top to bottom).\label{fig:spinod}}
\end{figure}

For the second example, we also use a piecewise constant approximation ($k=0$) with a Cartesian discretization of the domain $\Omega = (0,1)^2$, where $h=1.95\cdot 10^{-3}$. The interface parameter is taken to be very small $\gamma = 5\cdot 10^{-3}$, the time step is $\tau = 1\cdot 10^{-5}$ and $\mathrm{Pe}=50$. The initial condition for the order-parameter is taken to be a random value between $-1$ and $1$ inside a circular partition of the Cartesian mesh and $-1$ outside. The velocity field is given by
\begin{equation*}
\vec{u}(\vec{x}) := \frac12\left(1+\tanh (80-200\|(x_1-0.5,x_2-0.5) \|_{2})\right) \cdot \begin{pmatrix}2x_2-1\\ 1-2x_1 \end{pmatrix}, \quad \forall \vec{x}\in\Omega.
\end{equation*}
See Figure \ref{fig:spinod} for the numerical result. The method is robust with respect to $\gamma$ and is also well-suited to approach the thin high-gradient area of the order-parameter.


\subsection{Effect of the P\'eclet number\label{subsec:Pe}}

\begin{figure}[t]\center
\begin{minipage}[t]{0.99\columnwidth}\center
\includegraphics[width=3.2cm]{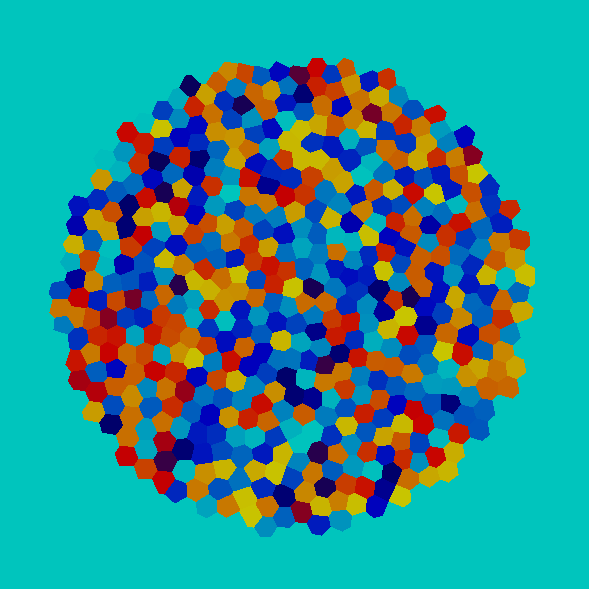}\hspace*{0.5cm}
\includegraphics[width=3.2cm]{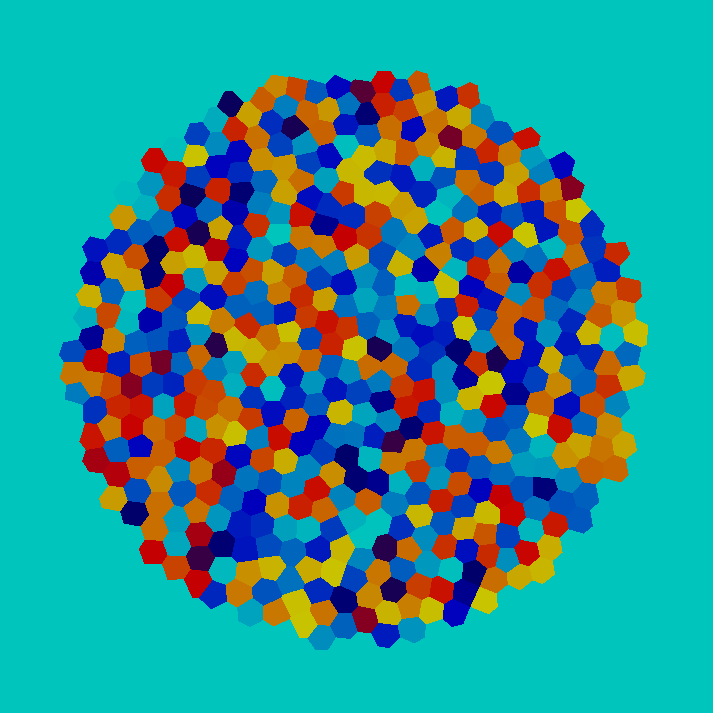}\hspace*{0.5cm}
\includegraphics[width=3.2cm]{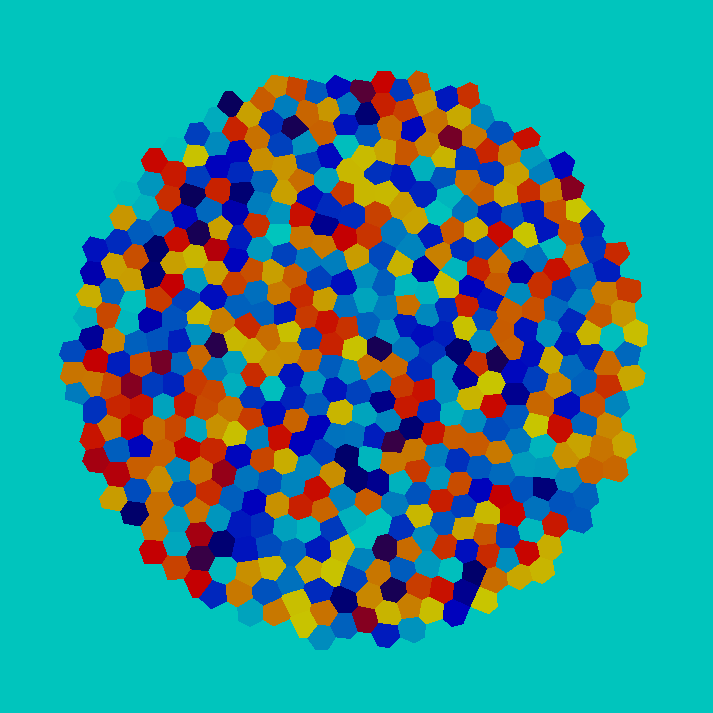}\vfill
\includegraphics[width=3.2cm]{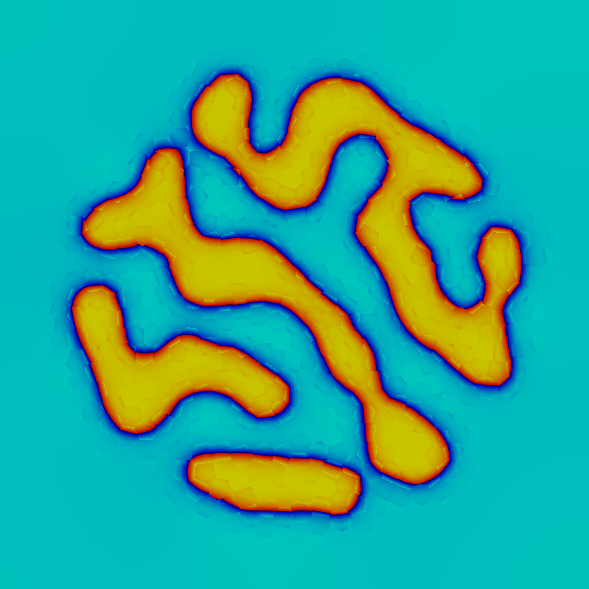}\hspace*{0.5cm}
\includegraphics[width=3.2cm]{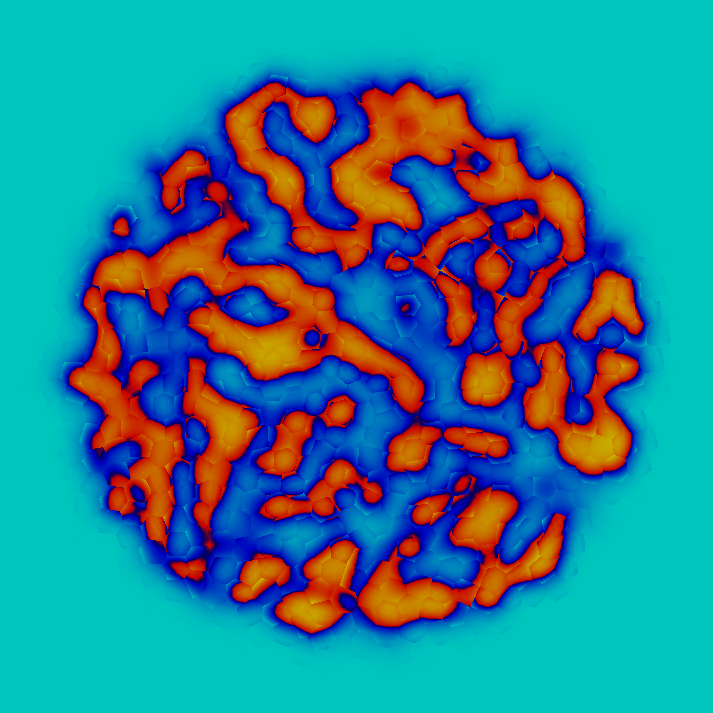}\hspace*{0.5cm}
\includegraphics[width=3.2cm]{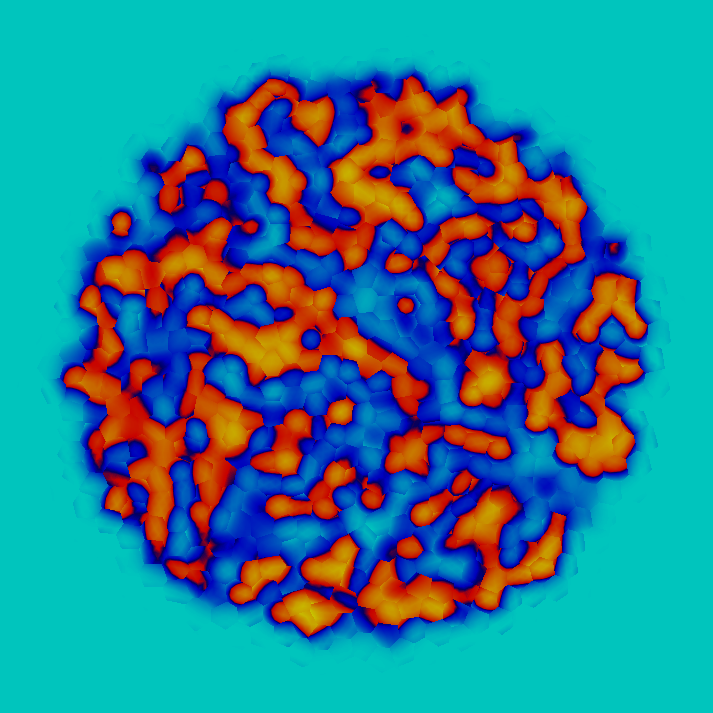}\vfill
\includegraphics[width=3.2cm]{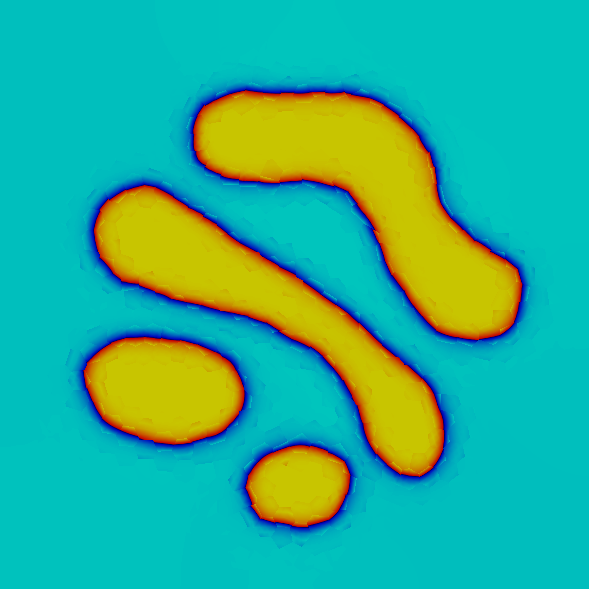}\hspace*{0.5cm}
\includegraphics[width=3.2cm]{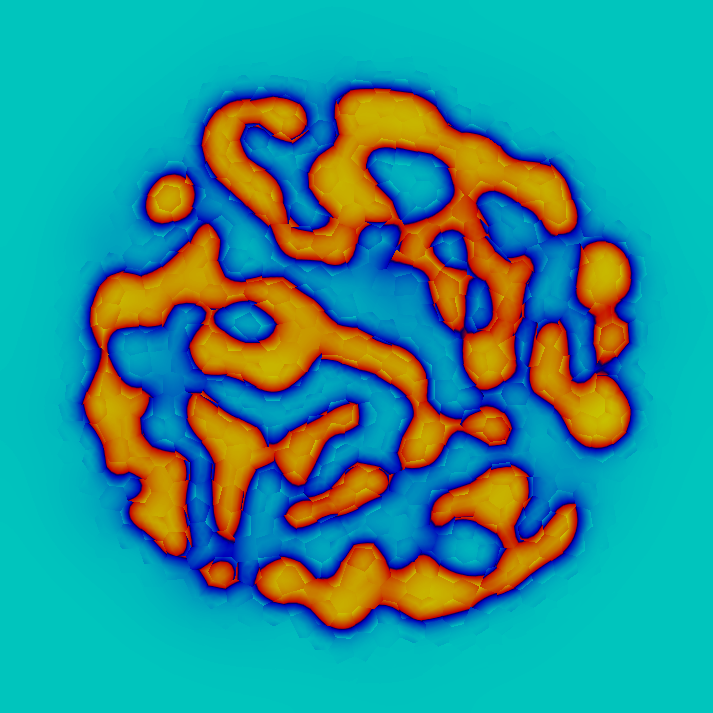}\hspace*{0.5cm}
\includegraphics[width=3.2cm]{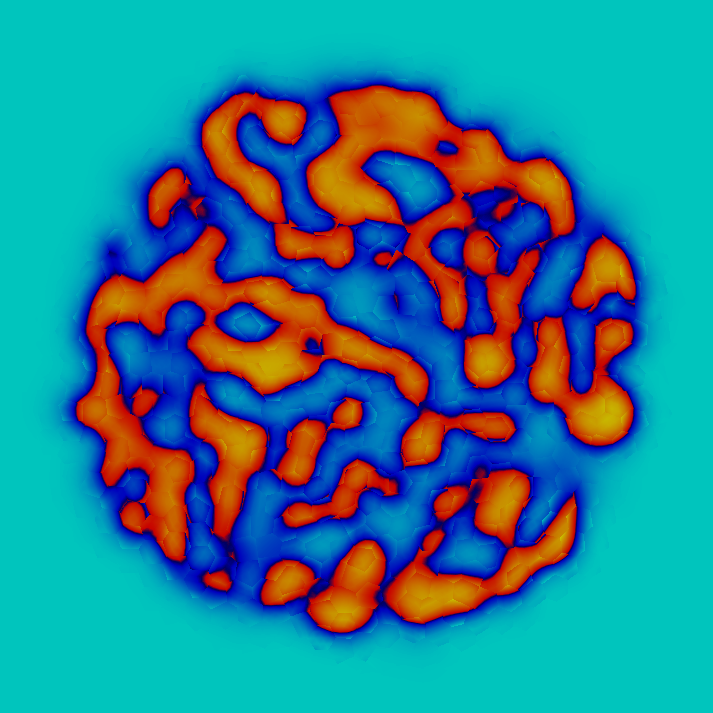}\vfill
\includegraphics[width=3.2cm]{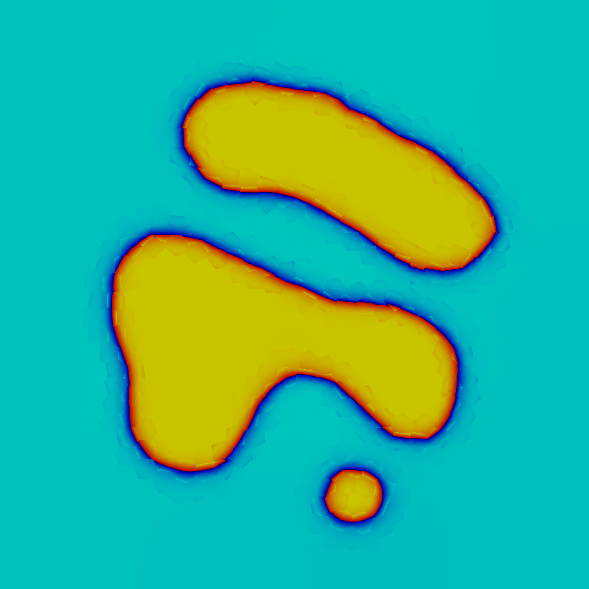}\hspace*{0.5cm}
\includegraphics[width=3.2cm]{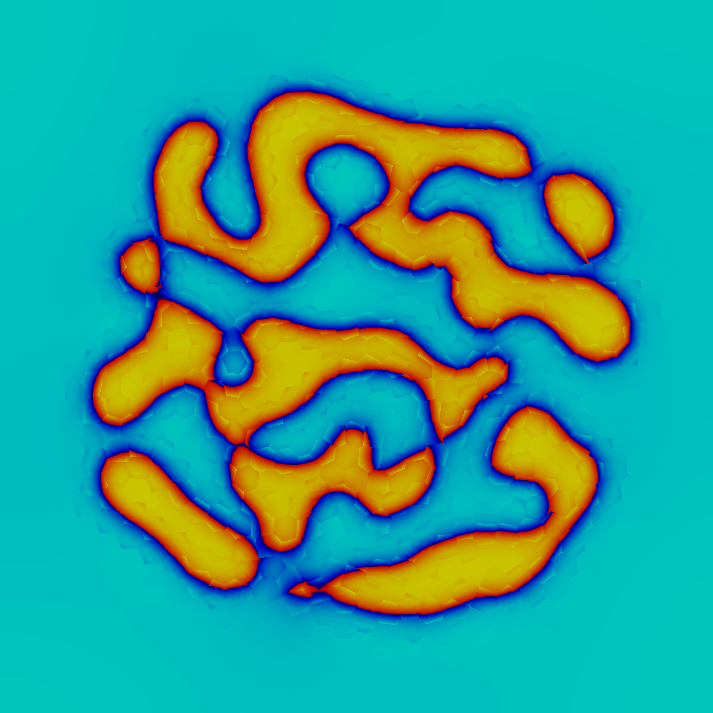}\hspace*{0.5cm}
\includegraphics[width=3.2cm]{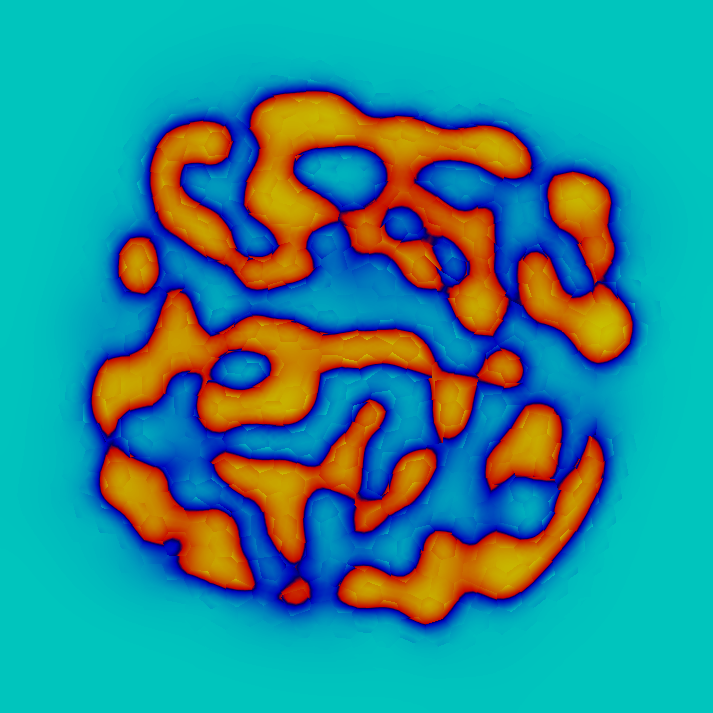}\vfill
\includegraphics[width=3.2cm]{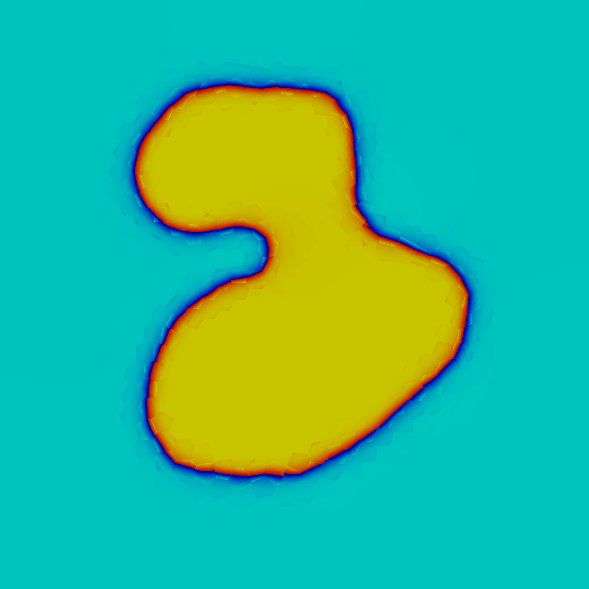}\hspace*{0.5cm}
\includegraphics[width=3.2cm]{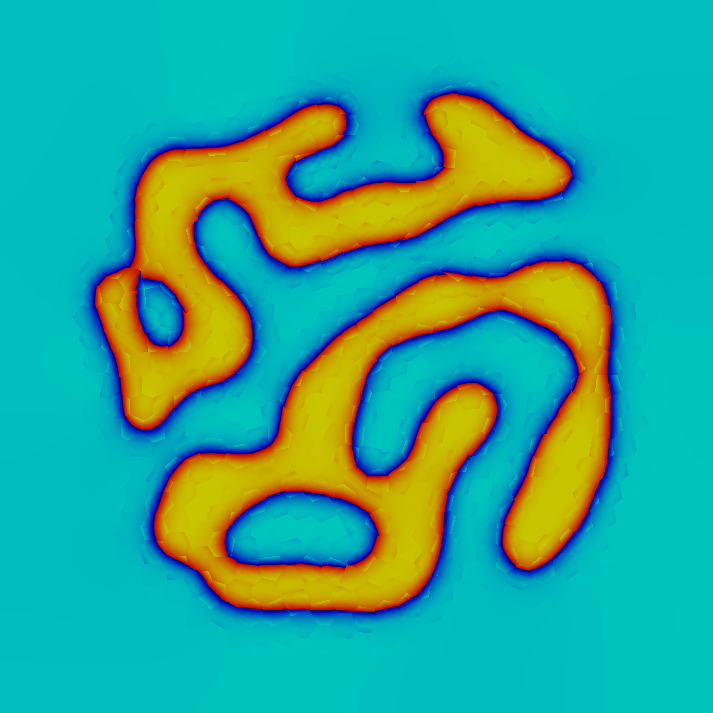}\hspace*{0.5cm}
\includegraphics[width=3.2cm]{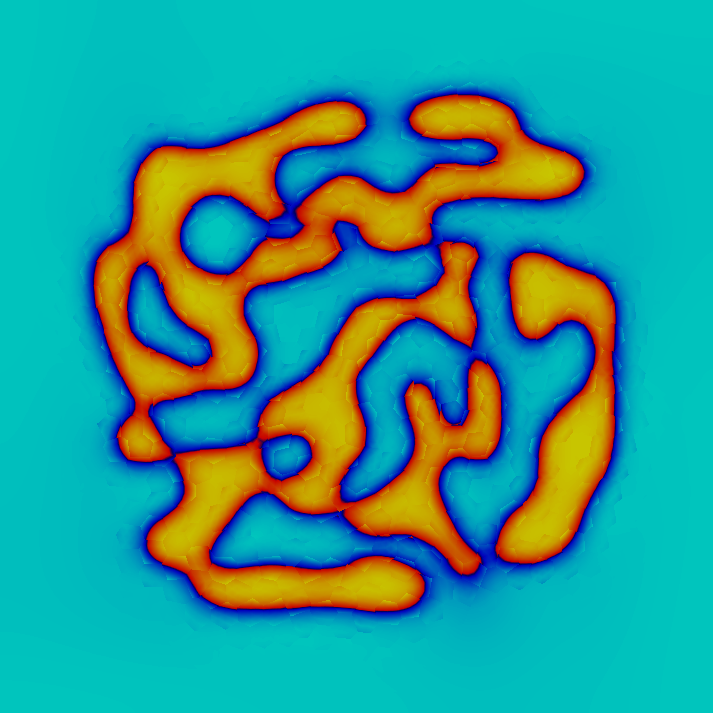}\vfill
\includegraphics[width=3.2cm]{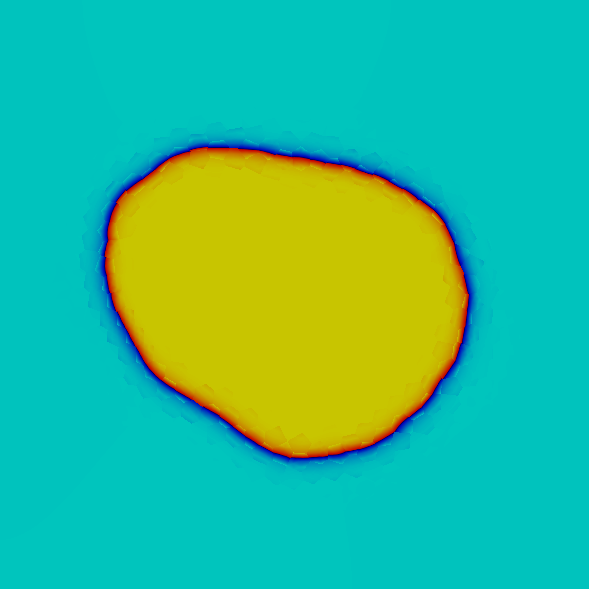}\hspace*{0.5cm}
\includegraphics[width=3.2cm]{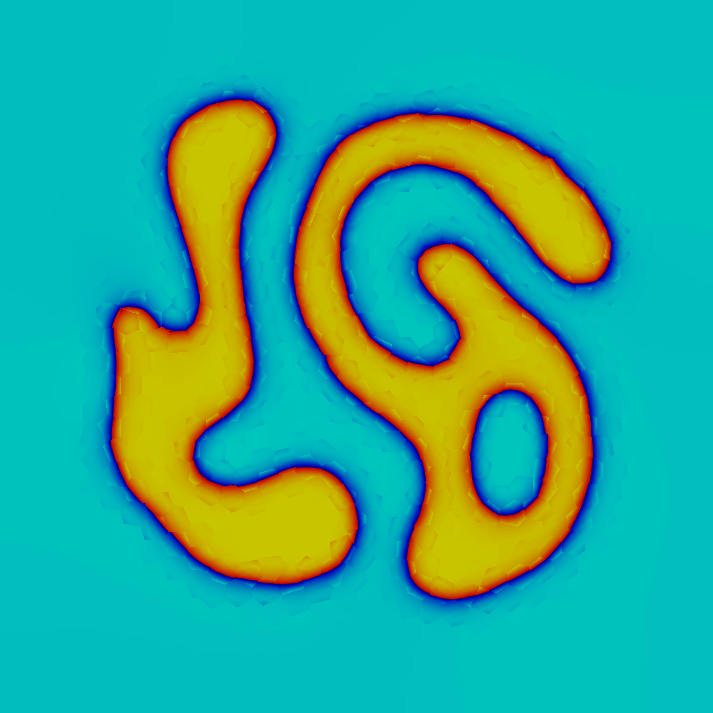}\hspace*{0.5cm}
\includegraphics[width=3.2cm]{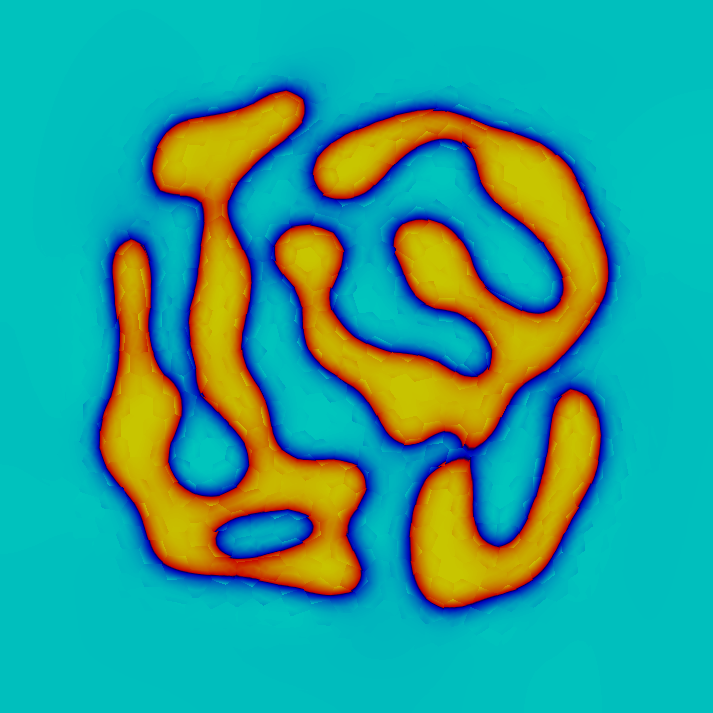}
\end{minipage}
\caption{Comparison at the same time between evolution of solutions with different P\'eclet number (top to bottom). Left: ${\rm Pe}=1$, middle: ${\rm Pe}=50$, right: ${\rm Pe}=200$. Displayed times are $t=0, 1\cdot 10^{-2}, 6\cdot 10^{-2}, 2\cdot 10^{-1}, 5\cdot 10^{-1}, 1$.\label{fig:comparison}}
\end{figure}

The P\'eclet number is the ratio of the contributions to mass transport by convection to those by diffusion: when ${\rm Pe}$ is greater than one, the effects of convection exceed those of diffusion in determining the overall mass flux. In the last test case, we compare several time evolutions obtained with different values of the P\'eclet number (${\rm Pe}\in\{1,50,200\}$), starting from the same initial condition. We use a Voronoi discretization of the domain $\Omega = (0,1)^2$, where $h=9.09\cdot 10^{-3}$, and use piecewise linear approximation ($k=1$). We choose $\gamma = 1\cdot 10^{-2}$, $\tau = 1\cdot 10^{-4}$ and $t_F = 1$. The initial condition is given by a random value between $-1$ and $1$ inside a circular domain of the Voronoi mesh and $-1$ outside. The convective term is given by
\begin{equation*}
\vec{u}(\vec{x}) := \begin{pmatrix}\sin(\pi x_1) \cos(\pi x_2)\\ -\cos(\pi x_1)\sin(\pi x_2) \end{pmatrix}, \quad \forall \vec{x}\in\Omega.
\end{equation*}
Snapshots of the order parameter at several times are shown on Figure \ref{fig:comparison} for each value of the P\'eclet number. For each case, the method takes into account the value of $\mathrm{Pe}$ and appropriately models the evolution of the order parameter by prevailing advection to diffusion when ${\rm Pe} \gg 1$.

\bibliographystyle{spmpsci}
\bibliography{chho_fvca.bib}

\begin{thebibliography}{1}
\providecommand{\url}[1]{{#1}}
\providecommand{\urlprefix}{URL }
\expandafter\ifx\csname urlstyle\endcsname\relax
  \providecommand{\doi}[1]{DOI~\discretionary{}{}{}#1}\else
  \providecommand{\doi}{DOI~\discretionary{}{}{}\begingroup
  \urlstyle{rm}\Url}\fi

\bibitem{Cahn:61}
Cahn, J.W.: On spinoidal decomposition.
\newblock Acta Metall. Mater. \textbf{9}, 795--801 (1961)

\bibitem{Cahn.Hilliard:58}
Cahn, J.W., Hilliard, J.E.: Free energy of a nonuniform system, {I},
  interfacial free energy.
\newblock J. Chem. Phys. \textbf{28}, 258--267 (1958)

\bibitem{Chave.Di-Pietro.ea:16}
Chave, F., Di~Pietro, D.A., Marche, F., Pigeonneau, F.: A hybrid high-order
  method for the {Cahn--Hilliard} problem in mixed form.
\newblock SIAM J. Numer. Anal. \textbf{54}(3), 1873--1898 (2016).
\newblock \doi{10.1137/15M1041055}

\bibitem{Di-Pietro.Droniou.ea:14}
Di~Pietro, D.A., Droniou, J., Ern, A.: A discontinuous-skeletal method for
  advection-diffusion-reaction on general meshes.
\newblock SIAM J. Numer. Anal. \textbf{53}(5), 2135--2157 (2015).
\newblock \doi{10.1137/140993971}

\bibitem{Di-Pietro.Ern:12}
Di~Pietro, D.A., Ern, A.: Mathematical aspects of discontinuous {G}alerkin
  methods, \emph{Mathématiques \& Applications}, vol.~69.
\newblock Springer-Verlag, Berlin (2012)

\end{thebibliography}

\begin{acknowledgement}
The work of D. A. Di Pietro and F. Marche was partially supported by \textit{Agence Nationale de la Recherche} grant HHOMM (ref. ANR-15-CE40-0005).
\end{acknowledgement}
\end{document}